\documentclass{ceb-l}

\usepackage{amssymb}

\usepackage{graphicx}
\usepackage{color}

\usepackage[cmtip,all]{xy}
\usepackage{tikz}
\usepackage{epigraph}
\usepackage{mathrsfs}

\newtheorem{theorem}{Theorem}[section]

\theoremstyle{definition}
\newtheorem{definition}[theorem]{Definition}
\newtheorem{example}[theorem]{Example}

\theoremstyle{remark}

\numberwithin{equation}{section}

\newcommand{\Sn}{\mathfrak{S}_n}
\newcommand{\PFn}{PF _n}
\newcommand{\suchthat}{\;|\;}
\newcommand{\area}{\ensuremath{\operatorname{area}}}
\newcommand{\dinv}{\ensuremath{\operatorname{dinv}}}
\newcommand{\Dinv}{\ensuremath{\operatorname{Dinv}}}
\newcommand{\word}{\ensuremath{\operatorname{word}}}
\newcommand{\ides}{\ensuremath{\operatorname{ides}}}
\newcommand{\Qsym}{\ensuremath{\operatorname{QSym}}}
\newcommand{\Sym}{\ensuremath{\operatorname{Sym}}}
\newcommand{\Des}{\ensuremath{\operatorname{Des}}}
\newcommand{\Mult}{\ensuremath{\operatorname{Mult}}}
\newcommand{\Char}{\ensuremath{\operatorname{Char}}}
\newcommand{\touch}{\ensuremath{\operatorname{touch}}}

\newlength\cellsize 
\savebox2{%
\begin{picture}(15,15)
\put(0,0){\line(1,0){15}}
\put(0,0){\line(0,1){15}}
\put(15,0){\line(0,1){15}}
\put(0,15){\line(1,0){15}}
\end{picture}}
\newcommand\cellify[1]{\def\thearg{#1}\def\nothing{}%
\ifx\thearg\nothing
\vrule width0pt height\cellsize depth0pt\else
\hbox to 0pt{\usebox2\hss}\fi%
\vbox to 15\unitlength{
\vss
\hbox to 15\unitlength{\hss$#1$\hss}
\vss}}
\newcommand\tableau[1]{\vtop{\let\\=\cr
\setlength\baselineskip{-16000pt}
\setlength\lineskiplimit{16000pt}
\setlength\lineskip{0pt}
\halign{&\cellify{##}\cr#1\crcr}}}
\savebox3{%
\begin{picture}(15,15)
\put(0,0){\line(1,0){15}}
\put(0,0){\line(0,1){15}}
\put(15,0){\line(0,1){15}}
\put(0,15){\line(1,0){15}}
\end{picture}}
\newcommand\expath[1]{%
\hbox to 0pt{\usebox3\hss}%
\vbox to 15\unitlength{
\vss
\hbox to 15\unitlength{\hss$#1$\hss}
\vss}}
\newcommand\bas[1]{\omit \vbox to \cellsize{ \vss \hbox to \cellsize{\hss$#1$\hss} \vss}}

\begin{document}

\title{The shuffle conjecture}


\author{Stephanie van Willigenburg}
\address{Department of Mathematics, University of British Columbia, Vancouver, BC V6T 1Z2, Canada}
\email{steph@math.ubc.ca}
\thanks{The author was supported in part by the National Sciences and Engineering Research Council of Canada, and in part by funding from the Simons
Foundation and the Centre de Recherches Math\'{e}matiques, through the Simons-CRM
scholar-in-residence program.}



\date{}

\dedicatory{On the occasion of Adriano Garsia's 90th birthday}

\begin{abstract} Walks in the plane taking unit-length steps north and east from $(0,0)$ to $(n,n)$ never dropping below $y=x$ and parking cars subject to preferences are two intriguing ingredients in a formula conjectured in 2005, now famously known as the shuffle conjecture.

Here we describe the combinatorial tools needed to state the conjecture. We also give key parts and people in its history, including its eventual algebraic solution by Carlsson and Mellit, which was published in the Journal of the American Mathematical Society in 2018. Finally, we conclude with some remaining open problems.
\end{abstract}

\maketitle


\epigraph{They can see the topography ... the treetops, but we can see the parakeets.}{Adriano Garsia}

Often, in order to delve deep into the structure of an abstract mathematical construct, the treetops, we need to interpret it concretely with a combinatorial visualization, the parakeets. The shuffle conjecture, as we will see, is one such story. In this article we will integrate the motivation, history and mathematics of the shuffle conjecture as we proceed. Hence we will begin by recalling necessary concepts from combinatorics in Section~\ref{sec:combintro}, and from algebra in Section~\ref{sec:algintro}, in order to state the shuffle conjecture in Theorem~\ref{the:shuffle}. This recently proved conjecture is, in essence, a formula for encoding the graded dimensions of the symmetric group representation in the character of a particular vector space on which the symmetric group $\Sn$ acts. In Section~\ref{sec:shuffle} we also discuss some of the motivation and history of the shuffle conjecture, including its refinement known as the compositional shuffle conjecture whose algebraic resolution by Carlsson and Mellit, announced in 2015 \cite{CMarxiv} and published in 2018 \cite{CarlssonMellit}, excited the combinatorial community. We mention some of their proof ingredients in Section~\ref{sec:further}, where we also conclude with some future avenues. 

\section{The combinatorics of Dyck paths and parking functions}\label{sec:combintro} A crucial concept for the statement of the shuffle conjecture is that of  parking functions. Although originally studied by Pyke \cite{Pyke}, they were introduced as a model for parking $n$ cars subject to preferences by Konheim and Weiss who were studying data storage \cite[Section 6: A parking problem -- the case of the capricious wives]{KonheimWeiss}. Konheim and Weiss also proved that the number of parking functions involving $n$ cars is $(n+1)^{(n-1)}$. Since then these functions have arisen in a plethora of places from hyperplane arrangements \cite{StanleyShi} to chip-firing \cite{CoriRossin}. More details on parking functions can be found, for example, in the survey by Yan \cite{Yan}. Rather than using the original definition, in terms of drivers parking cars, we will instead use an equivalent definition introduced by Garsia, for example in his paper with Haiman \cite[p 227]{GarsiaHaiman}. However, before we do this, we need to define a Dyck path.

\begin{definition}[Dyck path]\label{def:dyck} A \emph{Dyck path} of order $n$ is a path in the $n \times n$ lattice from $(0,0)$ to $(n,n)$ that consists of $n$ unit-length north steps and $n$ unit-length east steps, which stays weakly above the line $y=x$.
\end{definition}

\begin{example}\label{ex:dyck} If we let $N$ denote a unit-length north step, and $E$ denote a unit-length east step, then the following path $NNNEEENNENEENNEE$ from $(0,0)$ in the bottom-left corner to $(8,8)$ in the top-right corner is a Dyck path of order 8.

\

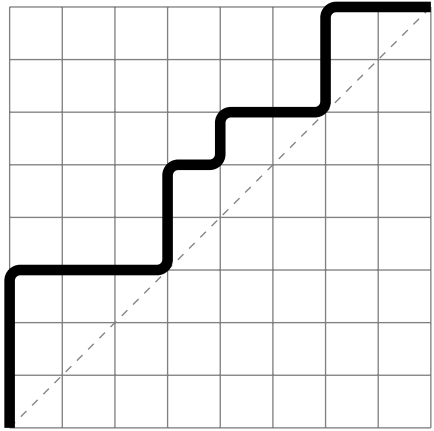
\begin{center}
    \begin{tikzpicture}[scale=0.7, every node/.style={font=\Large}]
        \useasboundingbox (-1,0) rectangle (8, 8);
        \draw[line width=0.5pt, color=gray] (0, 0) grid (8, 8);
        \draw[line width=0.5pt, dashed, color=gray] (0, 0) -- (8, 8);
        \draw[rounded corners=4, color=black, line width=4pt] (0, 0) -- (0, 1) -- (0, 2) -- (0, 3) -- (1, 3) -- (2, 3) -- (3, 3) -- (3, 4) -- (3, 5) -- (4, 5) -- (4, 6) -- (5, 6) -- (6, 6) -- (6, 7) -- (6, 8) -- (7, 8) -- (8, 8);
    \end{tikzpicture}
\end{center}
\end{example}

\begin{definition}[parking function]\label{def:PF} A \emph{parking function} of order $n$ is a Dyck path of order $n$ such that each north step has a label, called a \emph{car}, written in the square to its immediate right. The cars are $1, 2, \ldots , n$, each occurring exactly once, and cars in the same column increase from bottom to top. We denote the set of all parking functions of order $n$ by $\PFn$.
\end{definition}

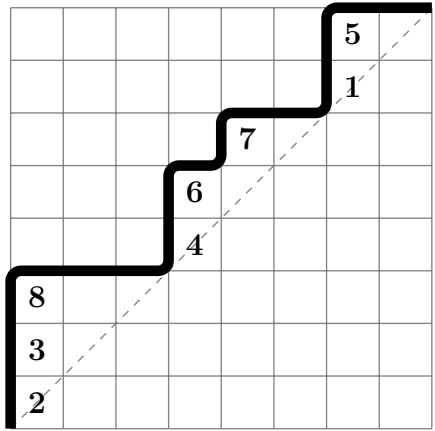
\begin{figure}[h]
    \centering
    \begin{tikzpicture}[scale=0.7, every node/.style={font=\Large}]
        \useasboundingbox (-1, 0) rectangle (8, 8);
        \draw[line width=0.5pt, color=gray] (0, 0) grid (8, 8);
        \draw[line width=0.5pt, dashed, color=gray] (0, 0) -- (8, 8);
        \draw[rounded corners=4, color=black, line width=4pt] (0, 0) -- (0, 1) -- (0, 2) -- (0, 3) -- (1, 3) -- (2, 3) -- (3, 3) -- (3, 4) -- (3, 5) -- (4, 5) -- (4, 6) -- (5, 6) -- (6, 6) -- (6, 7) -- (6, 8) -- (7, 8) -- (8, 8);
        \node[] at (0.5, 0.5) {$\mathbf{2}$};
        \node[] at (0.5, 1.5) {$\mathbf{3}$};
        \node[] at (0.5, 2.5) {$\mathbf{8}$};
        \node[] at (3.5, 3.5) {$\mathbf{4}$};
        \node[] at (3.5, 4.5) {$\mathbf{6}$};
        \node[] at (4.5, 5.5) {$\mathbf{7}$};
        \node[] at (6.5, 6.5) {$\mathbf{1}$};
        \node[] at (6.5, 7.5) {$\mathbf{5}$};
    \end{tikzpicture}
    \caption{A parking function of order 8}
    \label{fig:PF}
\end{figure}

\begin{example}\label{ex:PF}An example of a parking function, which we will use throughout this article, is given in Figure~\ref{fig:PF}.\end{example}

We now define three statistics on parking functions that will be useful later, the first of which is the area of a parking function and  depends only on its Dyck path.

\begin{definition}[\area]\label{def:area} If $\pi$ is a parking function, then its \emph{area} is the number of complete squares between the Dyck path of $\pi$ and $y=x$, denoted by $\area(\pi)$.
\end{definition}

\begin{example}\label{ex:area} If $\pi$ is the parking function from Figure~\ref{fig:PF}, then by counting the number of complete squares in each row contributing to the area, from bottom to top, we get
$$\area(\pi) = 0+1+2+0+1+1+0+1 = 6.$$
\end{example}

The second statistic is slightly more intricate than the area.

\begin{definition}[\dinv]\label{def:dinv}  Consider a parking function $\pi$, and a pair of cars $\{c_1, c_2\}$ in it.
\begin{itemize}
\item If the cars $c_1, c_2$ are in the \textbf{same} diagonal (that is, their squares are the same distance from $y=x$) with the larger car occurring further right, then $\{c_1, c_2\}$ is a \emph{primary diagonal inversion}. Let $\Dinv^{pri}(\pi)$ be the set of all such pairs.
\item If the cars $c_1, c_2$ are in \textbf{adjacent} diagonals with the larger car occurring in the higher diagonal (that is, its square is distance 1 further from $y=x$ than that of the smaller car) and further left, then $\{c_1, c_2\}$ is a \emph{secondary diagonal inversion}. Let $\Dinv^{sec}(\pi)$ be the set of all such pairs.
\end{itemize}
Then
$$\dinv(\pi) = | \Dinv^{pri}(\pi) | + | \Dinv^{sec}(\pi) |.$$
\end{definition}

\begin{example}\label{ex:dinv} If $\pi$ is from Figure~\ref{fig:PF}, then $\{3,7\}$ is a primary diagonal inversion, but $\{5,7\}$ is not, since the smaller car 5 occurs further right. Likewise $\{5,8\}$ is a secondary diagonal inversion, but $\{3,4\}$ is not, since the smaller car 3 occurs in the higher diagonal and further left. Note that $\{4,8\}$ is neither type of diagonal inversion since the cars are not in the same or adjacent diagonals.

Hence,
$$\Dinv^{pri}(\pi) = \{\{2,4\}, \{3,6\}, \{3,7\}, \{3,5\}, \{6,7\}\}$$
$$\Dinv^{sec}(\pi) = \{\{1,3\}, \{1,6\}, \{1,7\}, \{6,8\}, \{7,8\}, \{5,8\}\}$$so $$\dinv(\pi) = | \Dinv^{pri}(\pi) | + | \Dinv^{sec}(\pi) | = 5+6 = 11.$$
\end{example}

Our third statistic is a permutation associated to a parking function.

\begin{definition}[\word]\label{def:word} If $\pi$ is a parking function, then its \emph{word} is the permutation in one-line notation obtained by reading cars from the diagonal furthest from $y=x$ to the diagonal $y=x$, and within a diagonal reading from right to left. We denote this by $\word(\pi)$.
\end{definition}

\begin{example}\label{ex:word} If $\pi$ is from Figure~\ref{fig:PF}, then 
$$\word(\pi) = 85763142.$$
\end{example}

With our three statistics now defined, we end this section by recalling the $i$-descent set of a permutation, in our case specialized to the word of a parking function.

\begin{definition}[\ides]\label{def:ides} If $\pi$ is a parking function, then its \emph{$i$-descent set} is
$$\ides(\pi) = \{ i \suchthat i+1 \mbox{ is left of $i$ in $\word(\pi)$}\}.$$
\end{definition}

\begin{example}\label{ex:ides} If $\pi$ is from Figure~\ref{fig:PF} with $\word(\pi) = 85763142$ from Example~\ref{ex:word}, then 
$$\ides(\pi) = \{2,4,6,7\}.$$
\end{example}

\section{The algebras of quasisymmetric and symmetric functions}\label{sec:algintro} We now start to turn our attention to the algebraic ingredients needed to state the shuffle conjecture after first recalling the notions of compositions and partitions. 

A \emph{composition} $\alpha$ of $n$, denoted by $\alpha \vDash n$, is a list of positive integers $\alpha = \alpha _1\alpha _2 \cdots \alpha _{\ell(\alpha)}$ such that $\sum _{i=1} ^{\ell(\alpha)} \alpha _i = n$. We call the $\alpha _i$ the \emph{parts} of $\alpha$, call $n$ the \emph{size} of $\alpha$ and call $\ell(\alpha)$ the \emph{length} of $\alpha$. If, furthermore, $\alpha _1\geq\alpha _2\geq \cdots \geq\alpha _{\ell(\alpha)}$, then we say that $\alpha$ is a \emph{partition} of $n$, and denote this by $\alpha \vdash n$. For example, $332$ is both a composition and partition, with size 8 and length 3.

Now we focus on defining the algebra of quasisymmetric functions, before using them to define the algebra of symmetric functions.

The algebra of quasisymmetric functions, $\Qsym$, is a subalgebra of $\mathbb{C}[[z_1, z_2, \ldots ]]$, meaning that $\Qsym$ is a vector space, over $\mathbb{C}$, of formal power series in the variables $z_1, z_2, \ldots$, in which we can also multiply the elements together. A basis for it is given by the set of all fundamental quasisymmetric functions that we now define in the variables $Z=\{z_1, z_2, \ldots\}$, indexed by $n$ and subsets of $[n-1] = \{1,2, \ldots ,n-1\}$.

\begin{definition}[fundamental quasisymmetric function]\label{def:F} Let $S=\{ s_1, s_2, \ldots , s_{|S|}\}\subseteq [n-1]$. Then the \emph{fundamental quasisymmetric function} $F_{n,S}$ is defined to be
$$F_{n,S} = \sum z_{i_1}z_{i_2}\cdots z_{i_n}$$where the sum is over all $n$-tuples $(i_1, i_2, \ldots , i_n)$ satisfying
$$i_1\leq i_2 \leq \cdots \leq i_n \mbox{ and } i_j<i_{j+1} \mbox{ if } j\in S.$$
\end{definition}

\begin{example}\label{ex:F} We have that
$$F_{3, \{1\}} = z_1z_2^2+z_1z_3^2+ z_2z_3^2+\cdots +z_1z_2z_3+z_1z_2z_4+\cdots$$whereas
$$F_{3, \{2\}} = z_1^2z_2+z_1^2z_3+z_2^2z_3+\cdots +z_1z_2z_3+z_1z_2z_4+\cdots .$$
\end{example}

Quasisymmetric functions were first mentioned implicitly in Stanley's thesis,  with regard to $P$-partitions, published in 1972 \cite{Stanleythesis}, and then Gessel developed and published much of the classical theory explicitly in 1984 \cite{Gessel}. Since then they have arisen in a variety of areas, for example, from probability \cite{HershHsiao} to category theory \cite{AguiarBergeronSottile}. However, our interest lies in a special case of a result from Gessel's original paper \cite[Theorem 3]{Gessel}. For this we first need to define Young diagrams and Young tableaux.

Given a partition $\lambda = \lambda _1\lambda _2 \cdots \lambda _{\ell(\lambda)}\vdash n$, we define its \emph{Young diagram}, also denoted by $\lambda$, to be the array of left-justified boxes with $\lambda _i$ boxes in row $i$ from the top. Given the Young diagram, $\lambda$, a \emph{standard Young tableau (SYT)} of \emph{shape} $\lambda$, $T$, is a filling of the $n$ boxes of $\lambda$ with $1, 2, \ldots , n$ each appearing exactly once such that the entries in the rows  increase when read from left to right, and the entries in each column increase when read from top to bottom. We denote the set of all SYTs of shape $\lambda$ by $SYT(\lambda)$.

\begin{example}\label{ex:SYT} We have that $T = \tableau{1&3&4&5\\2&6&8\\7}$ is an SYT of shape $431\vdash 8$.
\end{example}

Given an SYT, $T$, of shape $\lambda \vdash n$, we define its \emph{descent set} to be
$$\Des (T) = \{ i \suchthat i+1 \mbox{ is in the same column or left of } i\} \subseteq [n-1].$$

\begin{example}\label{ex:des} If $T$ is from Example~\ref{ex:SYT}, then
$$\Des (T) = \{1,5,6\} \subseteq [7].$$
\end{example}

We can now define the algebra of symmetric functions, $\Sym$, which is a subalgebra of $\Qsym$. This algebra is so called because its elements are invariant under any permutation of its variables, and a basis for it is the set of all Schur functions that we now define as a special case of \cite[Theorem 3]{Gessel}.

\begin{definition}[Schur function]\label{def:schur} Let $\lambda \vdash n$. Then the \emph{Schur function} $s_\lambda$ is defined to be
$$s_\lambda = \sum _{T\in SYT(\lambda)} F_{n, \Des(T)}.$$
\end{definition}

\begin{example}\label{ex:schur} We have that $s_{21} = F_{3, \{1\}} + F_{3, \{2\}}$ from the SYTs below.
$$\tableau{1&3\\2}\qquad \tableau{1&2\\3}$$
\end{example}

The Schur functions are not the only basis of $\Sym$. Another basis that will be vital to our story is the basis consisting of all elementary symmetric functions: We define the \emph{$i$-th elementary symmetric function} $e_i$ to be
$$e_i = s _{1^i}$$where $1^i$ is the partition consisting of $i$ parts equal to 1. Then if $\lambda = \lambda _1 \lambda _2 \cdots \lambda _{\ell(\lambda)}\vdash n$ we define the \emph{elementary symmetric function} $e_\lambda$ to be
$$e_\lambda = e_{\lambda _1}e_{\lambda _2}\cdots e_{\lambda _{\ell(\lambda)}}.$$

Symmetric functions date back to Girard \cite{Girard}  in 1629, although Schur functions are much younger, dating to an 1815 paper of Cauchy \cite{Cauchy}. The Schur functions were named after Schur who in 1901 proved that they were characters of the irreducible polynomial representations of the general linear group \cite{Schur}, while standard Young tableaux were defined by Young in his 1928 publication \cite[p 258]{Young}. Substantial historical notes on this subject can be found in Stanley's second volume on enumerative combinatorics \cite[Chapter 7]{ECII},  which is also an excellent resource for symmetric functions and some related representation theory, as is the book by Sagan \cite{Sagan}.

\section{The space of diagonal harmonics and the shuffle conjecture}\label{sec:shuffle} With our essential combinatorial and algebraic notations now defined, we can begin to work towards our statement of the shuffle conjecture, which is about the vector space of diagonal harmonics. However, before we do that, let us define our desired space.

\begin{definition}[space of diagonal harmonics]\label{def:DHn}
Let $X_n = \{x_1, x_2, \ldots ,x_n\}$ and $Y_n = \{y_1, y_2, \ldots ,y_n\}$. Then the \emph{space of diagonal harmonics}, $DH_n$, is the vector space of polynomials in these variables, $f(X_n, Y_n)$, which satisfy
\begin{equation}\label{eq:DHn}
\partial ^a _{x_1}\partial ^b _{y_1} f(X_n, Y_n) + \partial ^a _{x_2}\partial ^b _{y_2} f(X_n, Y_n) + \cdots + \partial ^a _{x_n}\partial ^b _{y_n} f(X_n, Y_n) = 0
\end{equation}for all $a,b\geq 0$ and $a+b>0$. That is, 
$$DH_n = \{ f(X_n, Y_n) \in \mathbb{C} [X_n, Y_n] \suchthat \sum _{i=1}^{n} \partial ^a _{x_i}\partial ^b _{y_i} f(X_n, Y_n) = 0, \forall a,b\geq 0, a+b>0\}.$$\end{definition}

\begin{example}\label{ex:DH2} $DH_2$ consists of all polynomials $f(X_2, Y_2) = f(x_1, x_2, y_1, y_2)$ such that
$$\begin{array}{lcccll}
a+b=1&\mbox{gives}&\partial _{x_1} f(X_2, Y_2) + \partial _{x_2} f(X_2, Y_2)=0&\mbox{when} &a=1&b=0\\
&&\partial _{y_1} f(X_2, Y_2) + \partial _{y_2} f(X_2, Y_2)=0&&a=0&b=1\\
a+b=2&\mbox{gives}&\partial _{x_1}^2 f(X_2, Y_2) + \partial _{x_2}^2 f(X_2, Y_2)=0&\mbox{when} &a=2&b=0\\
&&\partial _{x_1}\partial _{y_1} f(X_2, Y_2) + \partial _{x_2}\partial _{y_2} f(X_2, Y_2)=0&&a=1&b=1\\
&&\partial _{y_1}^2 f(X_2, Y_2) + \partial _{y_2}^2 f(X_2, Y_2)=0&&a=0&b=2
\end{array}$$etc, and we can check that the solution set has basis $\{1, x_1-x_2, y_1-y_2\}$.
\end{example}

The symmetric group, $\Sn$, acts naturally on $DH_n$ by the \emph{diagonal action} that permutes the $X_n$ and $Y_n$ variables simultaneously. Namely, given $\sigma \in \Sn$ and $f(X_n, Y_n) = f(x_1, x_2, \ldots, x_n, y_1, y_2, \ldots , y_n)$ we have that
$$\sigma f(x_1, x_2, \ldots, x_n, y_1, y_2, \ldots , y_n) = f(x_{\sigma(1)}, x_{\sigma(2)}, \ldots, x_{\sigma(n)}, y_{\sigma(1)}, y_{\sigma(2)}, \ldots , y_{\sigma(n)}).$$By Equation~\eqref{eq:DHn} we see that if $f(X_n, Y_n) \in DH_n$, then $\sigma f(X_n, Y_n) \in DH_n$. Furthermore, if we let $DH_n^{c,d}$ be the subspace of $DH_n$ whose elements have total degree $c$ in the variables $x_1, x_2, \ldots , x_n$, and total degree $d$  in the variables $y_1, y_2, \ldots , y_n$, then if $f(X_n, Y_n) \in DH_n ^{c,d}$, then $\sigma f(X_n, Y_n) \in DH_n ^{c,d}$. This enables us to define the \emph{bigraded Frobenius characteristic} of $DH_n$ to be
\begin{equation}\label{eq:DHchar}DH_n[Z;q,t] = \sum _{c,d\geq 0} t^cq^d \sum _{\lambda\vdash n} s_\lambda \Mult(\chi ^\lambda, \Char DH_n ^{c,d})\end{equation}where, as before, $s_\lambda$ is a Schur function in the variables $Z=\{z_1, z_2, \ldots\}$ and $\Mult(\chi ^\lambda, \Char DH_n ^{c,d})$ is the multiplicity of the irreducible character of $\Sn$,  $\chi ^\lambda$, in the character of $DH_n ^{c,d}$  under the diagonal action of $\Sn$, $\Char DH_n ^{c,d}$. 

\subsection{The shuffle conjecture}\label{subsec:shuffle} The shuffle conjecture is a combinatorial formula for computing $DH_n[Z;q,t]$ in Equation~\eqref{eq:DHchar}, but before we give it and do an example we will briefly recount a skeletal history of what motivated it. More details on this fascinating story can be found in the excellent state-of-the-art survey article by Hicks \cite{Hicks}, and the illuminating texts by  Bergeron \cite{Bergeronbook}  and Haglund \cite{Haglundbook}.

In 1988 Kadell looked for \cite{Kadell}  and then Macdonald  found \cite{Macdonaldpoly}  a generalization of Schur functions, with additional parameters $q,t$, $P_\lambda[Z;q,t]$ where $\lambda \vdash n$.  This  generalization specialized to the Schur functions at $q=t$, and to other well-known functions such as the elementary symmetric functions and Hall-Littlewood functions,   which were likewise recovered by setting $q$ and $t$ to various values. {These functions were then transformed by Garsia and Haiman \cite[p 194]{GarsiaHaiman}, thus creating modified Macdonald polynomials $\widetilde{H} _\lambda [Z; q,t]$, which they hoped to prove were a positive linear combination of Schur functions. Proving this would imply the Macdonald positivity conjecture} dating from Macdonald's original work in 1988, which conjectured that Macdonald polynomials were a positive linear combination of Schur functions. In order to prove their conjecture, they defined vector spaces $\mathcal{H} _\lambda$ \cite{GarsiaHaimanPNAS}, now known as Garsia-Haiman modules, and conjectured that the bigraded Frobenius characteristic of $\mathcal{H} _\lambda$ was $\widetilde{H} _\lambda [Z; q,t]$. Moreover, they conjectured \cite[Conjecture 1]{GarsiaHaimanPNAS} that {irrespective} of $\lambda$ we have the following.
$$\dim (\mathcal{H} _\lambda) = n !$$This conjecture became known as the $n!$ conjecture, and both it and the Macdonald positivity conjecture were eventually proved by Haiman \cite[Theorem 3.2]{HaimanJAMS}.

At the same time Garsia and Haiman were studying  $DH_n$, which contains all the $\mathcal{H} _\lambda$ for $\lambda \vdash n$ as subspaces, and conjectured a formula for $DH_n[Z;q,t]$ in terms of the $\widetilde{H} _\lambda [Z; q,t]$.  Bergeron and Garsia noted that this formula was almost identical to the formula for the elementary symmetric functions $e_n$ in terms of $\widetilde{H} _\lambda [Z; q,t]$. More precisely, if the coefficient of $\widetilde{H} _\lambda [Z; q,t]$ in $e_n$ was $\mathscr{C}_\lambda$, then its conjectured coefficient in $DH_n[Z;q,t]$ was
$$t^{n(\lambda)}q^{n(\lambda ')}\mathscr{C}_\lambda$$where if $\lambda = \lambda _1 \lambda _2 \cdots \lambda _{\ell(\lambda)}$, then $n(\lambda) =\sum _{i=1} ^{\ell(\lambda)}\lambda _i(i-1)$ and $\lambda ' = \lambda ' _1 \lambda ' _2 \cdots \lambda ' _{\ell(\lambda ')}$ is the \emph{transpose} of $\lambda$, which is the partition created from $\lambda$ by setting
$$\lambda ' _i = \mbox{ number of parts of } \lambda \mbox{ that are }\geq i.$$For example, if $\lambda = 211$, then $\lambda ' = 31$. This inspired Bergeron and Garsia to officially define the \emph{nabla operator} in the paper  \cite[Equation (4.11)]{BergeronGarsia}  as follows.
$$\nabla \widetilde{H} _\lambda [Z; q,t] = t^{n(\lambda)}q^{n(\lambda ')}  \widetilde{H} _\lambda [Z; q,t]$$Hence when Haiman, using algebraic geometry, proved the conjectured formula for $DH_n[Z;q,t]$ \cite[Theorem 3.2]{Haiman}  this automatically yielded that \cite[Proposition 3.5]{Haiman}
\begin{equation}\label{eq:dequalsnabla} DH_n[Z;q,t] = \nabla e_n
\end{equation}since from above
$$e_n = \sum _{\lambda \vdash n} \mathscr{C}_\lambda \widetilde{H} _\lambda [Z; q,t]$$and now it was proved that $$DH_n[Z;q,t] = \sum _{\lambda \vdash n} t^{n(\lambda)}q^{n(\lambda ')}\mathscr{C}_\lambda \widetilde{H} _\lambda [Z; q,t].$$Haiman had also proved \cite[Proposition 3.6]{Haiman} that
$$\dim (DH_n) = (n+1) ^{(n-1)}.$$This supported the search for a collection of $(n+1) ^{(n-1)}$ objects, such as all parking functions of order $n$, along with statistics on them, in order to find a formula to compute $\nabla e_n$ more easily. The shuffle conjecture of Haglund, Haiman, Loehr, Remmel and Ulyanov \cite[Conjecture 3.1.2]{HHLRU} proposed such a formula, which we give now. This conjecture was proved recently, as a consequence of proving a refinement of it called the compositional shuffle conjecture,  by Carlsson and Mellit \cite[Theorem 7.5]{CarlssonMellit}.  However, many still refer to it as the shuffle conjecture, and hence we will too.

\begin{theorem}[the shuffle conjecture]\label{the:shuffle}
$$\nabla e_n = \sum _{\pi \in PF_n} t^{\area (\pi)} q^{\dinv(\pi)} F_{n, \ides(\pi)}$$
\end{theorem}

\begin{example}\label{ex:shuffle} Let us compute $n=2$. In order to compute $\nabla e_2$, we first need to calculate the elements of $PF_2$ that are as follows. 

\begin{center}
    \begin{tikzpicture}[scale=0.7, every node/.style={font=\Large}]
        \draw[line width=0.5pt, color=gray] (0, 0) grid (2, 2);
        \draw[line width=0.5pt, color=gray, dashed] (0, 0) -- (2, 2);
        \draw[rounded corners=4, color=black, line width=4pt] (0, 0) -- (0, 1) -- (0, 2) -- (1, 2) -- (2, 2);
        \node[] at (0.5, 0.5) {$\mathbf{1}$};
        \node[] at (0.5, 1.5) {$\mathbf{2}$};
        \node[] at (1, -1) {\large$\pi^{(1)}$};
    \end{tikzpicture}
    \qquad
    \begin{tikzpicture}[scale=0.7, every node/.style={font=\Large}]
        \draw[line width=0.5pt, color=gray] (0, 0) grid (2, 2);
        \draw[line width=0.5pt, color=gray, dashed] (0, 0) -- (2, 2);
        \draw[rounded corners=4, color=black, line width=4pt] (0, 0) -- (0, 1) -- (1, 1) -- (1, 2) -- (2, 2);
        \node[] at (0.5, 0.5) {$\mathbf{2}$};
        \node[] at (1.5, 1.5) {$\mathbf{1}$};
        \node[] at (1, -1) {\large$\pi^{(2)}$};
    \end{tikzpicture}
    \qquad
    \begin{tikzpicture}[scale=0.7, every node/.style={font=\Large}]
        \draw[line width=0.5pt, color=gray] (0, 0) grid (2, 2);
        \draw[line width=0.5pt, color=gray, dashed] (0, 0) -- (2, 2);
        \draw[rounded corners=4, color=black, line width=4pt] (0, 0) -- (0, 1) -- (1, 1) -- (1, 2) -- (2, 2);
        \node[] at (0.5, 0.5) {$\mathbf{1}$};
        \node[] at (1.5, 1.5) {$\mathbf{2}$};
        \node[] at (1, -1) {\large$\pi^{(3)}$};
    \end{tikzpicture}
\end{center}

They have 
$$\begin{array}{lccclcc}
\area(\pi ^{(1)})&=&1&\ &\dinv(\pi ^{(1)})&=&0\\
\area(\pi ^{(2)})&=&0&\ &\dinv(\pi ^{(2)})&=&0\\
\area(\pi ^{(3)})&=&0&\ &\dinv(\pi ^{(3)})&=&1\\
&&&&&&\\
\word(\pi ^{(1)})&=&21&\ &\ides(\pi ^{(1)})&=&\{1\}\\
\word(\pi ^{(2)})&=&12&\ &\ides(\pi ^{(2)})&=&\emptyset\\
\word(\pi ^{(3)})&=&21&\ &\ides(\pi ^{(3)})&=&\{1\}\\
\end{array}$$and hence
$$\nabla e_2 = tF_{2, \{1\}} + F_{2,\emptyset} + q F_{2,\{1\}} = F_{2,\emptyset} + (q+t) F_{2,\{1\}}.$$By Equation~\eqref{eq:dequalsnabla} and the definition of $DH_n[Z;q,t]$ in Equation~\eqref{eq:DHchar} we know that $\nabla e_2$ can be written as a positive linear combination of Schur functions. Using Definition~\ref{def:schur} we have that
$$s_2 = F_{2,\emptyset}\quad \mbox{and}\quad s_{11} = F_{2,\{1\}} $$from the respective SYTs $$\tableau{1&2} \quad \mbox{ and }\quad \tableau{1\\2}$$and hence
$$\nabla e_2 = s_2 + (q+t) s_{11}.$$

It is still an open problem to find a formula for $\nabla e_n$ that is a manifestly positive linear combination of Schur functions. 
\end{example}

We conclude this subsection with an indication of why the shuffle conjecture was so called. The name arose because  the coefficient of the monomial $z_1^{\lambda_1} z_2^{\lambda_2}\cdots z_{\ell(\lambda)}^{\lambda_{\ell(\lambda)}}$ in $\nabla e_n$ is equal to \cite[Corollary 3.3.1]{HHLRU} 
$$\sum  t^{\area(\pi)}q^{\dinv(\pi)}$$where the sum is over all ${\pi \in PF_n}$ such that $\word(\pi)$ is a \emph{shuffle} of the lists
$$[1, 2, \ldots, \lambda _1],\ [\lambda _1 +1, \lambda _1 +2, \ldots ,\lambda _1+\lambda_2], \ldots ,\ [m+ 1, m+2, \ldots , n]$$where $m= \sum _{i=1} ^{\ell(\lambda) -1}\lambda _i$, that is, within $\word(\pi)$ the numbers within each list appear in order when $\word(\pi)$ is read from left to right.

\begin{example}\label{ex:wordshuffle} Given the lists $[\mathit{1},\mathit{2}]$ and $[\mathbf{3},\mathbf{4}]$ note that 
$\mathit{1}\mathbf{3}\mathbf{4}\mathit{2}$ is a shuffle of the lists, but $\mathit{1}\mathbf{4}\mathbf{3}\mathit{2}$ is not since $\mathbf{3}$ and $\mathbf{4}$ are not in order.
\end{example}

\subsection{The compositional shuffle conjecture}\label{subsec:compshuffle} The conjecture that Carlsson and Mellit proved was not the shuffle conjecture from the previous subsection, but rather a refinement of it known as the compositional shuffle conjecture. This refinement by Haglund, Morse and Zabrocki \cite[Conjecture 4.5]{HMZ} centred around further symmetric functions $C_\alpha$, where $\alpha \vDash n$, that satisfy
$$e_n = \sum _{\alpha \vDash n} C_\alpha$$so that
\begin{equation}\label{eq:easC}\nabla e_n = \sum _{\alpha \vDash n} \nabla  C_\alpha\end{equation}
and involved a fourth statistic on parking functions, that of a touch composition.

\begin{definition}[\touch]\label{def:touch} If $\pi$ is a parking function of order $n$, then note the set of row numbers from bottom to top where there is a car in a square on the diagonal $y=x$
$$\{i_1=1, i_2, \ldots , i_k\}.$$Then the \emph{touch composition} is
$$\touch(\pi) = (i_2-i_1)(i_3-i_2)\cdots(n+1-i_k).$$
\end{definition}

\begin{example}\label{ex:touch} If $\pi$ is from Figure~\ref{fig:PF}, then the set of row numbers where there is a car in a square on $y=x$ is $\{1,4,7\}$ and hence
$$\touch(\pi) = 332.$$
\end{example}

We can now state the compositional shuffle conjecture of Haglund, Morse and Zabrocki \cite[Conjecture 4.5]{HMZ}, which was proved by Carlsson and Mellit \cite[Theorem 7.5]{CarlssonMellit}.

\begin{theorem}[the compositional shuffle conjecture] \label{the:compshuff} Let $\alpha \vDash n$.
$$\nabla C_\alpha = \sum _{\pi \in PF_n \atop \touch(\pi)=\alpha} t^{\area (\pi)} q^{\dinv(\pi)} F_{n, \ides(\pi)}$$
\end{theorem}

Observe that proving this would immediately prove the shuffle conjecture since if we sum over all $\alpha \vDash n$, then the left-hand side would yield $\nabla e_n$ by Equation~\eqref{eq:easC} and the right-hand side would lose its touch composition restriction.

\begin{example}\label{ex:compshuff} Let us compute $n=2$. From Example~\ref{ex:shuffle} we have that the elements of $PF_2$ are again as follows. 

\begin{center}
    \begin{tikzpicture}[scale=0.7, every node/.style={font=\Large}]
        \draw[line width=0.5pt, color=gray] (0, 0) grid (2, 2);
        \draw[line width=0.5pt, color=gray, dashed] (0, 0) -- (2, 2);
        \draw[rounded corners=4, color=black, line width=4pt] (0, 0) -- (0, 1) -- (0, 2) -- (1, 2) -- (2, 2);
        \node[] at (0.5, 0.5) {$\mathbf{1}$};
        \node[] at (0.5, 1.5) {$\mathbf{2}$};
        \node[] at (1, -1) {\large$\pi^{(1)}$};
    \end{tikzpicture}
    \qquad
    \begin{tikzpicture}[scale=0.7, every node/.style={font=\Large}]
        \draw[line width=0.5pt, color=gray] (0, 0) grid (2, 2);
        \draw[line width=0.5pt, color=gray, dashed] (0, 0) -- (2, 2);
        \draw[rounded corners=4, color=black, line width=4pt] (0, 0) -- (0, 1) -- (1, 1) -- (1, 2) -- (2, 2);
        \node[] at (0.5, 0.5) {$\mathbf{2}$};
        \node[] at (1.5, 1.5) {$\mathbf{1}$};
        \node[] at (1, -1) {\large$\pi^{(2)}$};
    \end{tikzpicture}
    \qquad
    \begin{tikzpicture}[scale=0.7, every node/.style={font=\Large}]
        \draw[line width=0.5pt, color=gray] (0, 0) grid (2, 2);
        \draw[line width=0.5pt, color=gray, dashed] (0, 0) -- (2, 2);
        \draw[rounded corners=4, color=black, line width=4pt] (0, 0) -- (0, 1) -- (1, 1) -- (1, 2) -- (2, 2);
        \node[] at (0.5, 0.5) {$\mathbf{1}$};
        \node[] at (1.5, 1.5) {$\mathbf{2}$};
        \node[] at (1, -1) {\large$\pi^{(3)}$};
    \end{tikzpicture}
\end{center}

They have 
$$\begin{array}{lcc}
\touch(\pi ^{(1)})&=&2\\
\touch(\pi ^{(2)})&=&11\\
\touch(\pi ^{(3)})&=&11\\
\end{array}$$hence
\begin{align*}
\nabla C_2&=tF_{2, \{1\}}\\
\nabla C_{11}&=F_{2,\emptyset} + q F_{2,\{1\}}
\end{align*}
and from Example~\ref{ex:shuffle}
$$\nabla e_2 = tF_{2, \{1\}} + F_{2,\emptyset} + q F_{2,\{1\}} = \nabla C_2 + \nabla C_{11}.$$
\end{example}

\section{The proof and further directions}\label{sec:further} On 25 August 2015 Carlsson and Mellit posted an article on the arXiv \cite{CMarxiv} titled simply ``A proof of the shuffle conjecture'', in which they proved the compositional shuffle conjecture, which in turn proved the shuffle conjecture. In their proof they refined the compositional shuffle conjecture yet further and proved this refinement. 

They worked with the right-hand side of the compositional shuffle conjecture under what is known as the $\zeta$ map, which takes a parking function $\pi$ to a new Dyck path with cars placed in the squares along $y=x$ such that when the cars are read from right to left we obtain $\word(\pi)$.  This required them to develop an analogue of $\touch$ that they called $\touch '$. They also worked with the reverse ordering of cars, so that, for example, in a parking function the cars in the same column \emph{decrease} when read from bottom to top. The list of other ingredients that they were required to create is impressive and included a generalization of the double affine Hecke algebra; partial Dyck paths; numerous operators including raising and lowering operators involving Hecke algebra operators and plethysm, and a modification of Demazure-Lusztig operators; and a recurrence that their refinement satisfied. 

To give a further idea of the complexity of the proof, this proof was almost 30 pages in length. In order to make it more accessible to combinatorialists, at the request of Garsia, Haglund and Xin expanded the proof, and their resulting article \cite{HaglundXin} is 60 pages in length.

\subsection{Further directions}\label{subsec:further} Carlsson and Mellit's proof of the shuffle conjecture was published in the Journal of the American Mathematical Society in 2018 \cite{CarlssonMellit}, but there remain many related open problems, some of which we now conclude with.

\begin{enumerate}
\item \textbf{A Schur-positive formula for $\nabla e_n$} By Equation~\eqref{eq:dequalsnabla} and  Equation~\eqref{eq:DHchar} we know that when we express $\nabla e_n$ as a linear combination of Schur functions
$$\nabla e_n = \sum _{c,d\geq 0} t^cq^d \sum _{\lambda \vdash n} \mathscr{D}_\lambda s_\lambda$$we have that the coefficients $\mathscr{D}_\lambda$ must be nonnegative integers since they are counting multiplicities. It remains an open problem to find a \emph{combinatorial} formula for the $\mathscr{D}_\lambda$, namely a formula that would compute them directly as nonnegative integers by counting a set of objects.
\item \textbf{Nabla on other symmetric functions} While the search for a combinatorial formula for $\nabla e_n$ has now been concluded with the proof of the shuffle conjecture, it remains to  prove the formula of Loehr and Warrington \cite[Conjecture 2.1]{LoehrWarringtonnablas} for
$$\nabla s_\lambda$$as the formula would generalize the result for $\nabla e_n$ since $e_n = s_{1^n}$. However, a conjecture of Loehr and Warrington \cite[p 667]{LoehrWarrington}  for
$$\nabla p_n$$where $p_n$ is the \emph{$n$-th power sum symmetric function}
$$p_n = z_1^n+ z_2^n+\cdots$$was recently proved by Sergel \cite[Theorem 4.11]{Sergel}  who has also conjectured the existence of a formula  \cite[Conjecture 3.1]{Sergelconj}  for
$$\nabla m_\lambda$$where $m_\lambda$ is the \emph{monomial symmetric function}
$$m_\lambda = \sum z_{i_1}^{\lambda _1}z_{i_2}^{\lambda _2}\cdots z_{i_{\ell(\lambda)}}^{\lambda _{\ell(\lambda)}}$$for $\lambda = \lambda _1\lambda _2\cdots\lambda _{\ell(\lambda)} \vdash n$ and the indices and monomials are distinct.
\item \textbf{A formula for $q,t$-Kostka polynomials} The modified Macdonald polynomials $\widetilde{H}_\lambda$, $\lambda \vdash n$, can be expanded as a linear combination of Schur functions
$$\widetilde{H}_\lambda = \sum _{\mu \vdash n} \widetilde{K}_{\mu\lambda}(q,t) s_\mu$$where the $\widetilde{K}_{\mu\lambda}(q,t)$ are known as $q,t$-Kostka polynomials. It is still an open problem to find a combinatorial formula for them, although such formulas have been found for $\lambda = m1^{n-m}$ by Stembridge \cite[Theorem 2.1]{Stembridge}, and also for $\lambda = 2^m1^{n-2m}$  by Fishel \cite[Theorem 1.1]{Fishel}, and others. Assaf, furthermore, has a theorem that enables the unification of these two cases \cite[Theorem 18]{Assaf}.\end{enumerate}

\section{Acknowledgments}\label{sec:ack} The author would like to thank H\'{e}l\`{e}ne Barcelo, David Eisenbud, Adriano Garsia, Jim Haglund, Angela Hicks, Richard Stanley and Mike Zabrocki for many fascinating conversations, and the Centre de Recherches Math\'{e}matiques and the Laboratoire de Combinatoire et d'Informatique Math\'{e}matique where many of the conversations and much of the subsequent writing took place thanks to a Simons CRM Professorship. She is also grateful to Niall Christie, Samantha Dahlberg, Jim Haglund, Angela Hicks, Franco Saliola, Mike Zabrocki and the referee for feedback on an earlier draft of the manuscript, and especially to Angela Hicks who additionally generously shared her draft of the history of the shuffle conjecture. She is also very grateful to  Angela Hicks and Franco Saliola for the code that created the diagrams.

\bibliographystyle{natbib}

\end{document}